\documentclass[12pt,a4paper]{amsart}
\newfont{\sheaf}{eusm10 scaled\magstep1}
\usepackage{amssymb}
\usepackage{graphics}
\usepackage{color}

\newcommand{\Proof}{{\it Proof. }}
\newcommand{\QED}{{\hfill $Q.E.D.$}}

\newtheorem{teo}{Theorem}[section]
\newtheorem{df}[teo]{Definition}

\newtheorem{question}[teo]{Question}
\newtheorem{lem}[teo]{Lemma}

\newtheorem{ex}[teo]{Example}
\newtheorem{oss}[teo]{Remark}
\newtheorem{prop}[teo]{Proposition}
\newtheorem{conj}[teo]{Conjecture}

\newcommand\sB{{\mathcal B}}

\def\Bbb{\bf}

\def\C{{\Bbb C}}

\newcommand{\CC}{\ensuremath{\mathbb{C}}}
\newcommand{\RR}{\ensuremath{\mathbb{R}}}
\newcommand{\ZZ}{\ensuremath{\mathbb{Z}}}
\newcommand{\QQ}{\ensuremath{\mathbb{Q}}}

\newcommand{\NN}{\ensuremath{\mathbb{N}}}

\newcommand{\M}{\ensuremath{\mathbb{M}}}

\newcommand{\PP}{\ensuremath{\mathbb{P}}}

\newcommand{\ra}{\ensuremath{\rightarrow}}

\newcommand{\SSS}{\ensuremath{\mathcal{S}}}
\newcommand{\AAA}{\ensuremath{\mathcal{A}}}
\newcommand{\F}{\ensuremath{\mathbb{F}}}

\def\eea{\end{eqnarray*}}
\def\bea{\begin{eqnarray*}}

\newcommand\dual{\mathrel{\raise3pt\hbox{$\underline{\mathrm{\thinspace d
\thinspace}}$}}}

\newcommand\qe{\ifhmode\unskip\nobreak\fi\quad $\Box$}       

\def\BOX{\hfill\lower.5\baselineskip\hbox{$\Box$}}

\newcommand\Z{\Bbb Z}

\input xy
\xyoption{all}
\CompileMatrices

\begin{document}
%
%
%
%
\title[Beauville surfaces and group theory]{Chebycheff and Belyi polynomials,
dessins d'enfants, Beauville
   surfaces and group theory}


\author{I. Bauer}
\address{Mathematisches Institut\\
Universit\" at Bayreuth\\
95440 Bayreuth\\
Germany}
\email{Ingrid.Bauer@uni-bayreuth.de}

\author{F. Catanese}
\address{Lehrstuhl Mathematik VIII\\
Universit\" at Bayreuth\\
95440 Bayreuth\\
Germany}
\email{Fabrizio.Catanese@uni-bayreuth.de}

\author{F. Grunewald}

\address{Mathematisches Institut
der\\
Heinrich-Heine-Universit\"at D\"usseldorf\\
  Universit\"atsstr. 1\\
D-40225 D\"usseldorf\\
Germany}
\email{grunewald@math.uni-duesseldorf.de}

\thanks{
We would like to
thank Fabio Tonoli for helping us with the pictures.}

  \begin{abstract}
We start discussing the group of automorphisms of the field of complex
  numbers, and describe,
in the special case of polynomials with only two critical values,
Grothendieck's program
of 'Dessins d' enfants', aiming at giving  representations of the
absolute Galois group.
We describe Chebycheff and Belyi polynomials, and other explicit examples.
As an illustration, we briefly treat difference and Schur polynomials.
Then we concentrate on a higher dimensional analogue of the triangle curves,
namely, Beauville surfaces and varieties isogenous to a product.
We describe their moduli spaces, and show how the study of these varieties
leads to new interesting questions in the theory of finite (simple) groups.
  \end{abstract}

\maketitle

\section{Introduction}
This article is an extended version of the plenary talk given by the second
author at the CIMMA 2005 in Almeria. It is intended for a rather
general audience and   for this reason we start with explaining
very elementary and wellknown facts: we apologize to the expert
reader, who will hopefully be satisfied by the  advanced part.

We only give very selected proofs, and we often prefer to introduce
concepts via  examples in order to be selfcontained and
elementary.

The general main theme of this article is the interplay between algebra and
geometry. Sometimes algebra helps to understand the geometry of
certain objects and vice versa. The absolute Galois group, i.e., the group of
field automorphisms of the algebraic closure of $\QQ$, which is a
basic object of interest in algebraic number theory, is still very
mysterious and since Grothendieck's proposal of ``dessins
d'enfants'' people try to understand it via suitable actions on
classes of geometrical objects.

In this paper we address the simplest possible such action,
namely the action of the absolute Galois group
on the space of (normalized) polynomials with
exactly two critical values, and on some combinatorial objects
related to them.

At the basis of  such a correspondence lies
the so called "algebraization" of (non constant)
holomorphic functions $f : C \ra \PP^1$ from a compact Riemann surface
$C$ to the projective line:  the  Riemann existence
theorem  describes $f$ completely through its set of critical values
and the related monodromy homomorphism
(the combinatorial object mentioned above).

The class of polynomials with two critical values is far from being
understood, there are however two series of such polynomials,
the so called Chebycheff polynomials, and the Belyi polynomials.

While the latter have a very simple algebraic description but are not
so well known, the Chebycheff polynomials are ubiquitous in many fields
of mathematics (probability, harmonic analysis,..) and everybody
has encountered them some time (see e.g. \cite{rivlin}, \cite{rivlin2}).

We try to describe the Chebycheff polynomials from a geometric point of view,
explaining its geometrical ties with the function $cos (z)$, which
also has exactly two critical values (viz.,  $\pm 1$).
The group of symmetries for $ cos (z)$ is the infinite dihedral group
of transformations $ z  \mapsto \pm z + n, n \in \ZZ$, which
has as fundamental domain an euclidean triangle of type
$(2,2,\infty )$, i.e., with angles $ \frac\pi 2, \frac\pi 2, 0$.

Our complete mastering of the function  $ cos (z)$ allows us
to achieve a complete and satisfactory description of the polynomials
with finite dihedral symmetry, namely the Chebycheff polynomials.

This is a key idea for the program we address: the knowledge
of the uniformizing triangular functions in the hyperbolic case
(these are functions on the upper half plane, with group
of symmetries determined by a hyperbolic triangle with
angles $ \frac\pi { m_1}, \frac\pi { m_2}, \frac\pi { m_3}$, cf. 
\cite{caratheodory},part $7$, chapter
$2$) sheds light on the class
of the so called {\em triangle
   curves}, on which the absolute Galois group acts.

  Triangle curves are pairs $(C,f )$ of a curve $C$ and a holomorphic
Galois covering map $f : C
\ra \PP^1$ ramified exactly over $\{0,1,\infty\}$.

In fact,  being ramified only in three points, triangle curves
   are {\em rigid} (i.e., they do not have any non trivial
   deformations) hence they are defined over
   $\overline{\QQ}$.  Belyi's famous
theorem  says almost the converse: every curve defined over
   $\overline{\QQ}$ admits a holomorphic
  map $f : C
\ra \PP^1$ ramified exactly over $\{0,1,\infty\}$.
Belyi's theorem in turn aroused the interest of Grothendieck
to the possible application of this result to the  construction of interesting
representations of the absolute Galois group, on the set of the so called
"dessins d'enfants", which we  briefly describe in section $4$.

Our main purpose is to explain how Grothendieck's sets of "dessins d'enfants"
can be more conveniently replaced by the sets of fundamental groups
of certain higher dimensional varieties $X$ which admit an unramified
covering which is isomorphic to a product of curves (these are called
{\em varieties isogenous to a product}). The Galois group acts
transforming such a variety $X$ to another variety of the same type,
but  whose fundamental
group need not be isomorphic to the fundamental group of $X$.

The last chapter is thus devoted to an overview of
\cite{beauvillesurf}, where the so-called {\em Beauville surfaces} are
studied. Beauville  surfaces are rigid surfaces which admit a finite
unramified covering which is a product of two algebraic curves.

Giving a  Beauville surface is essentially equivalent to giving two triangle
curves with the same group $G$, and in such a way
  that $G$ acts freely on the product of the two
triangle curves.

Beauville surfaces (and their higher dimensional analogues)
   not only provide
a wide class of surfaces quite
manageable in order to test conjectures, but also show how close algebra
and geometry are. The ease with which one can handle these surfaces is
based on the fact  that they
  are determined by discrete combinatorial data. Therefore one
can translate existence problems, geometric properties etc. in a purely
group theoretic language.

It was very fascinating for us that, studying
Beauville surfaces, we found out that many geometric questions are
closely related to some classical problems and conjectures in
the theory of finite
groups.

\section{Field automorphisms}

Let $K$ be any field and let $\phi : K \ra K$ be an automorphism of $K$:
then, since $\phi (x) = \phi ( 1 \cdot x) =\phi (1) \cdot \phi (x)$, it follows
that $\phi (1) = 1$, therefore $\phi (n) = n$ for all $n \in \NN$. If
$K$ has characteristic $0$, i.e., $ n \neq 0$ for all $n$,
thus $\QQ \subset K$, then $\phi |_{\QQ}  = Id_{\QQ}$.

The following exercise, quite surprising for first year students,
asserts that
the real numbers have no automorphisms except the identity, and the
complex numbers have
'too many'.

\begin{lem}
1) $Aut (\RR) = \{ Id \}$;

2) $|Aut (\CC)| = 2^{2^{\aleph_0}}$.
\end{lem}

\Proof
1) $\phi (a^2) = \phi (a)^2$, thus $\phi$ of a square is a square, so
$\phi$ carries the set of squares $\RR_+$ to itself. In particular,
$\phi $ is increasing, and since it is equal to the identity on $\QQ$,
we get that $\phi$ is the identity.

2) This follows from the fact that, if $\sB$ and $\sB '$ are two
transcendency bases,
i.e., maximal subsets of elements satisfying no nontrivial polynomial
equations,
then $\sB$ and $\sB '$ have the same cardinality $2^{\aleph_0}$, and for
any bijection $\phi '$ between $\sB$ and $\sB '$ there exists an
automorphism $\phi$ such that $\phi |_{\sB} = \phi'$.
\QED
\begin{oss}
1) Part two of the previous lemma is essentially the theorem of Steinitz: two
algebraically closed fields
are isomorphic iff they have the same characteristic and the same
absolute transcendence degree.

2) In practice, the automorphisms of $\CC$ that we do understand are
only the continuous ones,
the identity and the complex conjugation $\sigma$ (i.e., $\sigma
(z) := \bar{z} = x - i y$).

3) The field of algebraic numbers $\overline{\QQ}$ is the set $\{ z
\in \CC | \exists P \in \QQ[x] \backslash \{0\},
s.t. \ P (z) = 0 \}$. The fact that $ Aut (\CC)$ is so large is due
to the fact that
the kernel of $Aut (\CC) \ra Aut (\overline{\QQ})$ is very large.

The group $Aut (\overline{\QQ})$ is called the {\em absolute Galois group}
and denoted by
$Gal(\overline{\QQ}/ \QQ)$.

Note that even if we have a presentation of a
group $G$, still our information about it might be quite scarce (even
the question:
``is the group nontrivial?'' is hard to answer), and the solution is
to have a {\em representation} of it, for instance, an action on a
set $M$ that can be
very well described.
\end{oss}

\begin{ex}
The dihedral group $D_n$ is best described through its action on $\CC$,
as the set of $2n$ transformations of the form
   $z \mapsto \zeta z $, where $\zeta^n = 1$, or
of the form $ z \mapsto \zeta \bar{z} $.
   We see immediately that the group acts
as a group of permutations of the regular
$n-$gon whose set of vertices is the set $\mu_n$ of
   $n-$th roots of unity $\mu_n : = \{\zeta |\zeta^n = 1 \}$.

The group is generated by complex conjugation $\sigma$, and by the rotation
$ r (z) = exp (\frac{2}{n} \pi i) z $, so its presentation
is $ < r, \sigma | r^n= 1, \sigma ^2 = 1, \sigma r \sigma  = r^{-1}
>$.
\end{ex}
We end the section by invoking
{\em Grothendieck's dream of dessins d' enfants (= children's drawings)},
which aims at finding  concrete representations for the Galois Group $Gal
(\overline{\QQ}/\QQ)$. We will try to explain some of the basic ideas
in the next sections, trying to be as elementary as possible.

One can invoke as a catchword {\em moduli theory}, in order to soon
impress the audience,
but it is possible to explain everything in a very simple way, which
is precisely what we shall do in the next sections, at least in some
concrete examples.


\section{Polynomials with rational critical values.}

Let $ P(z) \in \CC [z]$, $P(z) = \sum_{i =0}^n \ a_i z ^i$ be a polynomial of
degree $n$.
\begin{df}
1) $\zeta \in \CC$ is a {\em critical point} for $P$ if $ P ' ( \zeta ) = 0$;

2) $ w \in \CC$ is a {\em critical value} for $P$ if there is a  critical point
$\zeta$ for $P$ such that $w = P (\zeta)$.
Denote by $B_P$ the set of critical values
of $P$; $B_P$ is called the {\em branch set } of $P$.
\end{df}

\begin{oss}
1) $\phi \in Aut (\CC)$ acts on $ \CC [z]$, by $P (z) = \sum_{i =0}^n
  a_i z ^i
\mapsto  \phi (P) (z) : = \sum_{i =0}^n  \phi (a_i) z ^i$.

2) If $\zeta \in \CC$ is a  critical point for $P$, then, since
$\sum_{i =0}^n  i a_i \zeta ^{i-1} = 0$, $\phi (\zeta)$ is a critical
point for $\phi (P) $. Similarly, if $ w \in \CC$ is a  critical value for $P$,
then $\phi (w) = \phi (P (\zeta))= \phi (P) (\phi (\zeta))$ is a
critical value for $\phi (P)$.
\end{oss}

Since for each $\phi \in Aut (\CC)$ we have $\phi (\QQ) = \QQ$, and 
$\phi (\overline{\QQ}) =
\overline{\QQ}$, the class
of polynomials  $\{ P | B_P \subset \overline{\QQ} \}$, and
$\{ P | B_P \subset \QQ \}$ are invariant under the action
of $ Aut (\CC)$.

\begin{oss}
If $g : \CC \ra \CC $, $g(z) = az + b$ is an invertible affine 
transformation, i.e., $a \neq 0$, then two {\em right affine 
equivalent} polynomials
$P(z)$ and $P(g(z))$ are immediately
seen to have the same branch set. In order not to have infinitely
many polynomials with the same branch set, one considers the
{\em normalized polynomials} of degree $n$:
$$ P(z) = z^n + a_{n-2} z ^{n-2} + \dots   + a_0.$$

\noindent
Any polynomial is right affine equivalent to a normalized polynomial, 
and two normalized  polynomials
are right affine equivalent iff they are equivalent under the group 
$\mu_n$ of n-th roots of unity,
$P(z) \cong P(\zeta z)$ (i.e., $ a_i \mapsto \zeta^i a_i$).
\end{oss}

The usefulness of the above concept stems from the following

\begin{teo}\label{branch}
Let $P$ be a normalized polynomial: then
$P \in \overline{\QQ}[z]$ if and only if $B_P \subset \overline{\QQ}$.
\end{teo}
\Proof
It is clear that $P \in \overline{\QQ}[z]$ implies that $B_P \subset
\overline{\QQ}$.

Now, assume that $B_P \subset \overline{\QQ}$ and observe that,
by the argument of Steinitz, it suffices to show that $P \in
\CC[z]$
has only a finite number of transforms under the group $Aut(\CC)$.

Since a transform of a normalized polynomial is normalized, it suffices to show
that the right affine equivalence class of $P$ has only a finite number
of transforms under the group $ Aut(\CC)$.

Since the set $B_P  \subset \overline{\QQ}$ has only a finite number
of transforms,
it suffices to show that there is only a finite number of classes of
polynomials
having a fixed branch set $B$: but this follows from the well known
Riemann's existence theorem, which we shall now recall.
\QED

\begin{teo} (Riemann's existence theorem)
There is a natural bijection between:

1) Equivalence classes of holomorphic mappings  $f : C \ra \PP^1_{\CC}$,
   of degree $n$ and with branch set $B_f \subset B$,
(where $C$ is a compact Riemann surface, and $f : C \ra \PP^1_{\CC}$,
$f' : C' \ra \PP^1_{\CC}$ are said to be equivalent if there is a
biholomorphism
$ g : C' \ra C$ such that $ f' = f \circ g $).

2) Conjugacy classes of {\em monodromy} homomorphisms $ \mu : \pi_1
(\PP^1_{\CC} - B) \ra \SSS_n$
(here, $\SSS_n$ is the symmetric group in n letters, and $\mu \cong
\mu '$ iff there
is an element $\tau \in \SSS_n$ with  $\mu  (\gamma) = \tau \mu
'(\gamma) \tau^{-1}$, for all $\gamma$).

Moreover:

i) $C$ is connected if and only if the subgroup $ Im (\mu)$ acts transitively
on $\{1,2, \dots n \}$.

ii) $f$ is a polynomial if and only if $ \infty \in B$, the monodromy
at $ \infty$
is a cyclical permutation, and $g(C) =0$.
\end{teo}

For a proof we refer to \cite{miranda}, especially pp. 91,92.

\begin{oss}
1) Assume that $\infty \in B$, so write $B:= \{\infty , b_1, \dots b_d
\}$. Then
$\pi_1 (\PP^1_{\CC} - B)$ is a free group generated by
$\gamma_1, \dots \gamma_d$ and $\mu$ is completely determined by
the local monodromies $\tau_i : = \mu (\gamma_i)$.

2) To give then a polynomial
of degree $n$ with branch set $B_P = \{ b_1, \dots b_d \} $ it suffices to
give nontrivial permutations $\tau_1, \dots \tau_d$ such that
$$ \tau_1 \cdot  \dots  \cdot \tau_d  = (1,2, \cdots ,n),$$
and such that, if we write each $\tau_j$ as a product of disjoint cycles of
length $m_{ij}$, then $\sum_{i,j} m_{ij}= n-1$.

3) Riemann's existence theorem holds more generally also for maps of
infinite degree,
and it was generalized by Grauert and Remmert
(cf. \cite{grauertremmert}) to describe locally finite
holomorphic maps between normal complex spaces.

\end{oss}


\section{Polynomials with two critical values}

Observe that there is only one normalized polynomial of degree $n$
with only $0$ as
critical value, namely $P(z)=z^n$. In fact, there is only one choice for
$\tau_1$: $\tau_1   = (1,2, \cdots ,n).$

If we consider polynomials $P$ with only two critical values, we shall assume,
without loss of generality (and for historical reasons)
that $B_P =  \{ -1,1 \}$  or $B_P =  \{ 0,1 \}$.

We have to give two permutations $\tau_1, \tau_2$ such that
$\tau_1  \cdot \tau_2  = (1,2, \cdots ,n)$.
If  $\tau_1$ has a cycle decomposition of type $m_1, \dots, m_r$,
$\tau_2$ has a cycle decomposition of type $d_1, \dots ,d_s$,
then, counting the roots of the derivative of $P$, we find
that
$$
(**) \ \ \ \ \ \sum_j (m_j -1) + \sum_i (d_i - 1) = n-1.
$$

\begin{teo}
There is only one class of a polynomial $T_n$ of degree $n$,
called {\em Chebycheff polynomial of degree n}, such that

a) the critical values of $T_n$ are just $ \{ -1,1 \}$,

b) the derivative of $T_n$ has simple roots.

We have
$$
T_n (z) = \sum _{2r \leq n} (-1)^r z^{n-2r}
(1 - z^2)^r =
$$
$$
= cos ( n (arcos (z))  ) = \frac{1}{2} ( (z \pm \sqrt {z^2 -1} \ )^n
+ (z \pm \sqrt {z^2 -1} \ )^{-n}).
$$
\end{teo}

\begin{figure}[htbp]
\begin{center}
\scalebox{1}{\includegraphics{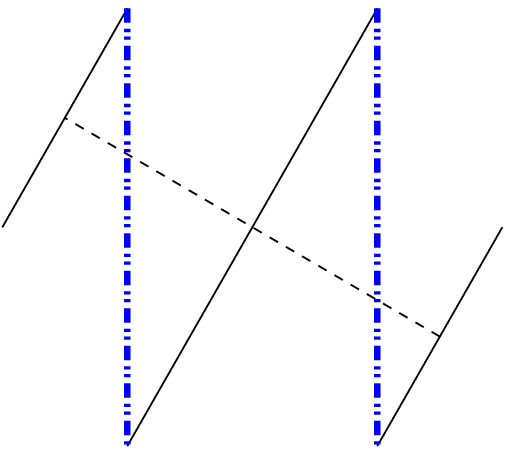}}
\qquad\qquad
\scalebox{1}{\includegraphics{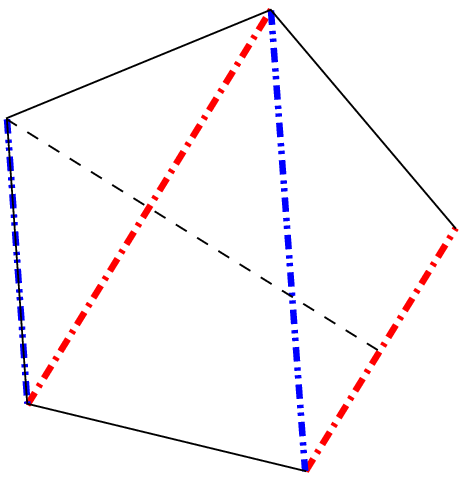}}
\end{center}
\caption{Chebycheff monodromies: even and odd}
\label{chebmonodromy}
\end{figure}

\Proof
By condition b), the permutations $\tau_1$, $\tau_2$ are a product
of disjoint transpositions, so $m_j = n_i = 2$ and formula $(**)$
says that we have exactly $n-1$ transpositions: one sees immediately
that there is only one combinatorial solution (cf. figure \ref{chebmonodromy}).

The monodromy image is in fact exactly the Dihedral group $D_n$:
this explains the above fomulas, because if we set $ u := T_n(z)$,
then there is a quadratic extension of $\CC (z)$, yielding the Galois closure
of $\CC (u) \subset \CC (z)$.

We have, if $w : = t^n$, $u = \frac{1}{2}
(w + w^{-1}) = \frac{w^2 + 1}{2w}$, and similarly
$z = \frac{t^2 + 1}{2t}$, thus $t^2 - 2 zt + 1 = 0$.

\begin{figure}[htbp]
\begin{center}
\input{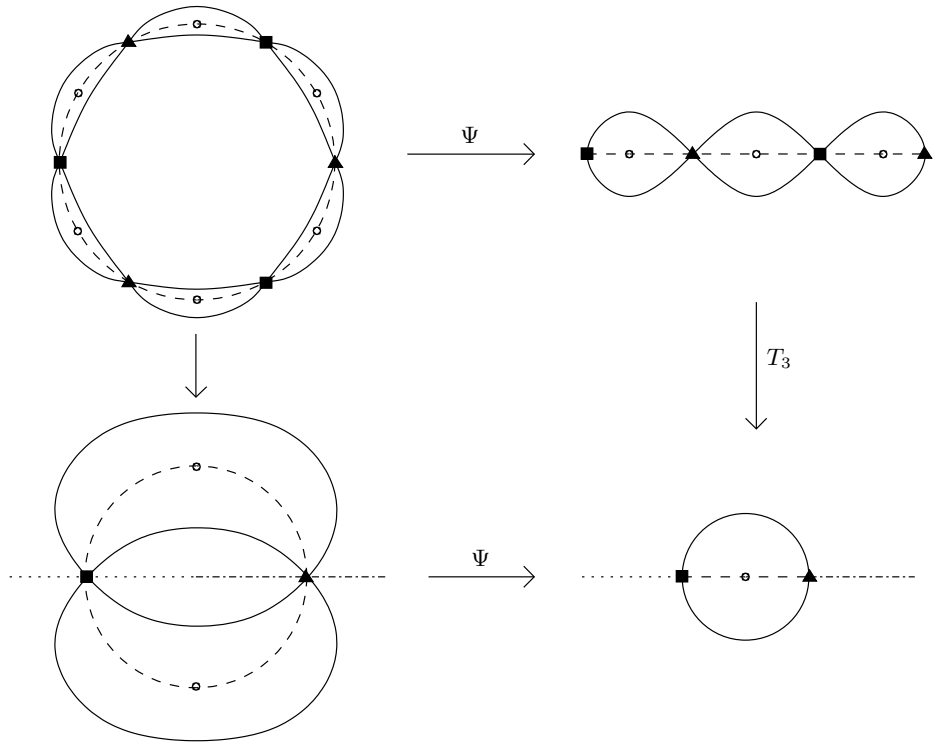}
\end{center}
\caption{Chebycheff basic diagram}
\label{basic}
\end{figure}

Therefore, we have the following basic diagram (cf. also figure \ref{basic}).
\begin{equation}\label{diagrtildestar}
\xymatrix{
\Psi: t\ar[r]\ar[d]&z = \frac{1}{2}(t + t^{-1}) = \frac{t^2+1}{2t}\ar_{T_n}[d]
\\
\Psi: t^n = w\ar[r]&u  = T_n(z) = \Psi(w).
}
\end{equation}

The relation with the function $cos (y)= \frac{1}{2} (e^{iy} +
e^{-iy})$ comes in
because the monodromy of $T_n$ factors through the mondromy of $cos$,
namely, the infinite Dihedral group
$$D_{\infty} = < a,b | a^2 = b^2 = 1> = AL (1, \ZZ): = \{ g | g(z) =
\pm z + c , \ c \in \ZZ\}.$$
\QED

Note that in the first picture of figure \ref{chebmonodromy} ($n = 
6$) we have two disjoint
transpositions corresponding to the reflection of the hexagon with
horizontal axis, and three disjoint transpositions corresponding to
the reflection with axis the dotted line.

In the second picture ($n=5$), we have two disjoint transpositions
corresponding to the reflection with horizontal axis and two other
disjoint transpositions corresponding to the reflection with axis the
dotted line.

Another very simple class of polynomials with two critical values is given
by the {\em Belyi polynomials}.

Take a rational number $q \in \QQ$ such that $0 < q < 1$: writing $q$ 
as a fraction $q =
\frac{m}{m+r}$, we get a polynomial
$$
P_{m,r} : = z^m (1-z)^r
(m+r)^{m+r}  m^{-m}  r^{-r}.
$$
The critical points of $ P_{m,r} $ are just $\{ 0,1, q = \frac{m}{m+r} \}$,
and its only critical values are $\{ 0,1 \}$. Observe that $P_{m,r} 
(0) = P_{m,r}(1) =0$,
while the coefficient of $P_{m,r}$ is chosen exactly in order that 
$P_{m,r} (q) = 1$.
Hence it follows that the monodromy $\tau'_1$  corresponding to
the critical value $1$ is a transposition, whereas
the monodromy $\tau'_0$  corresponding to
the critical value $0$ has a cycle decomposition of type $(m,r)$
(but be aware of the possibility  $m=1, r=1$).

\begin{df}
Let $P \in \CC [z]$ be a polynomial with critical values
$\{0,1\}$, and  observe
that $ 1 - P$ is also  a polynomial with critical values
$\{0,1\}$. We say that $ 1 - P$ is {\em extendedly}
equivalent to $P$, and we shall call the union of
the equivalence class of $P$ with the equivalence
class of  $ 1 - P$ the {\em extended equivalence class}.

\end{df}

In particular, for  degree $\leq 4$, Chebycheff polynomials
are extendedly equivalent to Belyi polynomials.

\begin{prop}\label{polynomials}

Consider the extended equivalence classes of polynomials with two 
critical values:
among them are the classes
of Chebycheff
and Belyi polynomials.

i) In degree $\leq 4$ there are no other classes.

ii)  In degree $5$ there is also
the class of the polynomial $\frac{3}{16} (z^5 - \frac{10}{3}z^3 + 5 
z + \frac{8}{3})$.

Observing that all the above polynomials are in $\QQ [z]$,
we have also that

iii) in degree $6$ we first find a polynomial $P$
which is not extendedly equivalent to any polynomial in $\RR [z]$,
and which more precisely is not extendedly equivalent
to its complex conjugate.

\end{prop}

\Proof
i) We leave as an exercise to the interested reader to find out
that all the possible (extended classes of) monodromies in degree $\leq 4$ are
either Chebycheff monodromies ($\iff$ $\tau_j$ for $ j=1,2$ is a product
of disjoint transpositions) or Belyi monodromies
($\iff$ $\tau_1$ is a transposition, and $\tau_2$ is a product of at 
most two cycles).

ii) Similarly one sees that in degree $5$ the only remaining case is the one
where each $\tau_j$ is a 3-cycle. Up to affine transformations we
may assume that the sum of the critical points equals $0$,
so that we get the polynomial
$$
Q (z) = 5 \int (z - a)^2 (z + a )^2 =
z^5 - \frac{10}{3} a^2 z^3 + 5 a^ 4 z +c.
$$
Imposing the condition $Q(-a) = 0$, $Q(a) = 1$ we obtain $c = 
\frac{1}{2}$, $a^5 = \frac{3}{16}$.

Thus we obtain a normalized polynomial which has coefficients which 
do not lie in $\QQ$.
However, if we set $ a=1$ and we multiply $ Q$ by $\frac3 {16}$, we
obtain a nonnormalized polynomial $ P \in \QQ [z]$ with
critical values $\{0,1\}$.

\begin{figure}[htbp]
\begin{center}
\scalebox{1}{\includegraphics{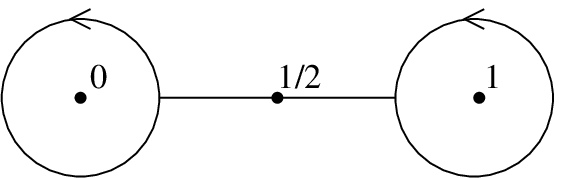}}
\end{center}
\caption{}
\label{dessin}
\end{figure}

To see iii)  is not difficult: it suffices to exhibit a polynomial $P$
whose monodromy is not conjugate to the monodromy
of $\bar{P}$. We choose now (cf. figure \ref{dessin}) $\frac12$
as base point, and basis of the fundamental group of
$\CC \setminus \{0,1\}$ such that complex conjugation
sends $\gamma_j \mapsto \gamma_j ^{-1}, \ j= 1,2$.

It follows that if $\tau_j , \ j= 1,2$ yield the monodromy
of $P$ , then $\tau_j ^{-1}, \ j= 1,2$ yield the monodromy
of $\bar{P}$. We choose $\tau_1 = (1,3, 6) (4,5), \tau_2 = (1,2) (3,5)$,
and observe that if there were a permutation $\alpha$ conjugating 
$\tau_i$ to $\tau_i^{-1}$
for $ i=1,2$,  then
$\alpha$ would leave invariant the sets $\{1,3,6\}$, $\{4,5\}$, 
$\{1,2\}$, $\{3,5\}$, hence their mutual
intersections.  This immediately
implies that $\alpha = id$, a contradiction.

Another example is given by setting $\tau_1 = (5,6) (1,2,3) , \tau_2 
=  (3,4,6)$.

\QED

Let $P \in \CC[z]$ be a polynomial with critical values
$\{0,1\}$. Then we know by theorem \ref{branch} that $P$ has
coefficients in $\overline{\QQ}$, and in fact then in some number
field $K$. We want to show this fact again, but  in a more explicit way
which will allow us to find effectively the field $K$.

Let us now introduce several affine algebraic sets which capture the
information contained in all the number fields arising this way.

Let $P(z):= z^n + a_{n-2}z^{n-2} + \ldots a_0$ be a normalized
polynomial with only critical values $\{0,1\}$. Once we choose the
types of the respective cycle decompositions $(m_1, \ldots, m_r)$ and $(n_1,
\ldots, n_s)$, we can write our polynomial $P$
in two ways, namely
$$
P(z) = \prod_{i=1}^r (z - \beta_i)^{m_i},
$$
$$
P(z)-1 = \prod_{k=1}^s (z - \gamma_k)^{n_k}.
$$
Since $P$ is normalized we have the equations $F_1 = \sum m_i \beta_i = 0$
and $F_2=\sum n_k \gamma_k = 0$.
We have: $m_1 + \ldots m_r = n_1 + \ldots n_s = n = degP$ and
therefore, since by $(**)$   $\sum_j (m_j-1) + \sum_i(n_i -1) = n-1$,
we get $r+s = n+1$.

Since we
have
$\prod_{i=1}^r (z -
\beta_i)^{m_i} = 1+ \prod_{k=1}^s (z - \gamma_k)^{n_k}$ comparing coefficients
we obtain further $n-1$ polynomial equations with integer coefficients in
the variables $\beta_i$, $\gamma_k$ which we denote by $F_3=0, \ldots,
F_{n+1}=0$.
Let
$${\bf V}\left(n,
\left[\begin{array}{ccc} m_1 & \ldots & m_r\\ n_1 & \ldots & 
n_s\end{array}\right]\right) $$
be the
algebraic set in affine $(n+1)$-space corresponding to the
equations $F_1=0, \ldots,F_{n+1}=0$. Mapping a point of this algebraic set to
the vector $(a_0,\ldots,a_{n-1})$ of coefficients of the corresponding
polynomial $P$ we obtain (by elimination of variables) an algebraic set
$${\bf C}\left(n,
\left[\begin{array}{ccc} m_1 & \ldots & m_r\\ n_1 & \ldots & 
n_s\end{array}\right]\right) $$
in affine $(n-1)$ space.
Both these algebraic sets are defined over $\QQ$. This
does not imply
that each of their points has coordinates in $\QQ$ (cf. proposition
\ref{polynomials}).
We encounter this
phenomenon also later, in the more complicated situation of Beauville
surfaces (cf. theorem \ref{nonreal}). We have
\begin{prop}
Each of the algebraic sets
$${\bf V}\left(n,
\left[\begin{array}{ccc} m_1 & \ldots & m_r\\ n_1 & \ldots &
     n_s\end{array}\right]\right),\quad
{\bf C}\left(n,
\left[\begin{array}{ccc} m_1 & \ldots & m_r\\ n_1 & \ldots &
     n_s\end{array}\right]\right) $$
is either empty or has dimension $0$.
\end{prop}
This follows from Riemann's existence theorem since this theorem shows that
there are only a finite
number of normalized polynomials having critical points $\{0,1\}$ and fixed
ramification types. We are going to describe this set of polynomials
  in the following
  case $n=6$, $r=4$, $s=3$. We shall give now information on the
algebraic sets $\bf C$ for all possible choices of the ramification indices.

\medskip
\centerline{{\bf I.} $(m_1,m_2,m_3,m_4)=(2,2,1,1)$, $(n_1,n_2,n_3)=(2,2,2)$:}

The algebraic set $\bf C$ is defined by the equations
$$a_0=-1,\quad a_1=0, \quad 4a_2=a_4^2, \quad a_3=0, \quad a_4^3=-54,$$
it is irreducible over $\QQ$ and consists of $3$ points over $\bar \QQ$
(in the Chebycheff class).

\medskip
\centerline{{\bf II.} $(m_1,m_2,m_3,m_4)=(2,2,1,1)$, $(n_1,n_2,n_3)=(3,2,1)$:}

The algebraic set $\bf C$ is defined by the equations
$$a_0= \frac{a_4^6}{8748} + \frac{a_4^ 3}{81} - \frac{2}{3},\ \
     a_1= - \frac{a_3a_4^4}{324} + \frac{a_3a_4}{3},\ \
     a_2= -\frac{25a_4^8}{314928} + \frac{2a_4^5}{729} + \frac{8a_4^2}{27},
$$
$$
     a_3^2 + \frac{5a_4^6}{2187} - \frac{4a_4^3}{81} - \frac{4}{3}=0,\ \
     a_4^9 - \frac{324a_4^6}{25} - \frac{17496a_4^3}{25} - \frac{314928}{25}=0,
$$
it is irreducible over $\QQ$ and consists of $18$ points over $\bar \QQ$.

\medskip
\centerline{{\bf III.} $(m_1,m_2,m_3,m_4)=(2,2,1,1)$, $(n_1,n_2,n_3)=(4,1,1)$:}

Here the algebraic set $\bf C$ falls into two irreducible components
${\bf C}_1,$ ${\bf C}_2$ already over $\QQ$. The variety ${\bf C}_1$
is defined by
$$a_0=-1,\quad a_1=0, \quad a_2=0, \quad a_3=0, \quad a_4^3=\frac{27}{4},$$
it is irreducible over $\QQ$ and consists of $3$ points over $\bar \QQ$.
The variety ${\bf C}_2$
is defined by
$$a_0 =- \frac{2527}{2500},\ \
a_1=  \frac{19a_3a_4}{150},\ \
         a_2= \frac{21a_4^2}{20},\ \
         a_3^2= - \frac{243}{125},\ \
         a_4^3 + \frac{27}{50}=0,
$$
it is irreducible over $\QQ$ and consists of $6$ points over $\bar \QQ$.

\medskip
\centerline{{\bf IV.} $(m_1,m_2,m_3,m_4)=(3,1,1,1)$, $(n_1,n_2,n_3)=(2,2,2)$:}

Here the algebraic set $\bf C$ falls into two irreducible components
${\bf C}_1,$ ${\bf C}_2$ already over $\QQ$. The components
${\bf C}_1,$ ${\bf C}_2$ are defined by
$$a_0=0,\quad a_1=0, \quad a_2=0, \quad a_3=\pm 2, \quad a_4=0.$$

\medskip
\centerline{{\bf V.} $(m_1,m_2,m_3,m_4)=(3,1,1,1)$, $(n_1,n_2,n_3)=(3,2,1)$:}

The algebraic set $\bf C$ is defined by the equations
$$
a_0=\frac{29a_4^3}{2700} -\frac{9}{50},\quad
     a_1=  \frac{7a_3a_4^4}{432} + \frac{3a_3a_4}{8},\quad
     a_2= - \frac{5a_4^5}{324} + \frac{5a_4^2}{12},$$
$$
a_3^2=  \frac{11a_4^3}{135} - \frac{2}{5},\quad
     a_4^6 - \frac{324a_4^3}{25} + \frac{2916}{25}=0,
$$
it is irreducible over $\QQ$ and consists of $12$ points over $\bar \QQ$.

\medskip
\centerline{{\bf VI.} $(m_1,m_2,m_3,m_4)=(3,1,1,1)$, $(n_1,n_2,n_3)=(4,1,1)$:}

The algebraic set $\bf C$ is defined by the equations
$$
a_0= - \frac{513}{625},\ \
     a_1= \frac{24a_3a_4}{25},\ \
     a_2=\frac{16a_4^2}{45},\ \
     a_3^2= \frac{32}{125},\ \
     a_4^3=- \frac{729}{100},
$$
it is irreducible over $\QQ$ and consists of $6$ points over $\bar \QQ$.
\bigskip

We return now to the general situation and we consider  the following set

$M_n : =  \{$right affine equivalence classes  of polynomials $P$ of
degree $n$
with critical
values $\{0,1\}\}$. As seen in thm. \ref{branch}, the Galois group
$Gal ( \overline{\QQ}/ \QQ )$ acts on $M_n$, and Grothendieck's expression:
``Dessins d' enfants'' refers to the graph $ P^{-1} ([0,1])$.

It has $n$ edges, corresponding to the inverse image of the open
interval $(0,1)$, and its set of vertices is {\em bipartite},
since we have two different types of vertices, say the black vertices
given by $ P^{-1} (\{0\})$, and
the blue vertices corresponding to $ P^{-1} (\{1\})$.

Moreover, around each vertex, we can give a cyclical counterclockwise
ordering of the edges incident in the vertex (cf. figure
\ref{dessin}).

It is easy to show that these data completely determine
the monodromy, which acts on the set $ P^{-1} (\{\frac{1}{2}\})$,
i.e., on the midpoints of the edges.
In fact, e.g. the ``black monodromy'' is given as follows: for each midpoint
of a segment:
go first to the black vertex, then turn to the right, until you reach the
next midpoint.

\begin{figure}[htbp]
\begin{center}
\scalebox{1}{\includegraphics{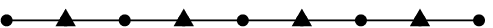}}
\end{center}
\caption{Dessins d'enfants: Chebycheff polynomial}
\label{dessincheb}
\end{figure}

\begin{figure}[htbp]
\begin{center}
\scalebox{1}{\includegraphics{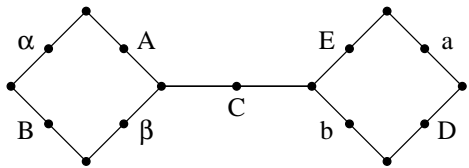}}
\end{center}
\caption{A rational function of degree $9$}
\label{dessindeg9}
\end{figure}

Figure \ref{dessincheb} gives the dessin d'enfant for a Chebycheff
polynomial, whereas figure \ref{dessindeg9} is the children 's drawing
of a rational function of degree $9$,
whose monodromy at $\infty$ is $(A,B,C,D,E)(a,b)(\alpha,\beta)$.

\begin{oss}
The action of the Galois group is still very mysterious, the only
obviously clear fact is that the conjugacy classes of
$\tau_1$, resp. $\tau_2$ (in other words, the multiplicities
of the critical points in $ P^{-1} (\{0\})$, resp. $ P^{-1} (\{1\})$)
remain unchanged.
\end{oss}

\section{ Riemann's existence theorem and difference
   polynomials}
In this section we briefly discuss some applications
  of Riemann's existence theorem, and in particular of Chebycheff
polynomials  to the problem of irreducibility
of difference
polynomials (cf. \cite{davenport}, \cite{schinzel}, \cite{fried70}, 
\cite{fried73}).

\begin{df}
Let $f$, $g \in \mathbb{C}[x]$ be two polynomials. Then a {\em
   difference polynomial} is a polynomial in $\mathbb{C}[x,y]$ of the
   form $f(x) - g(y)$.
\end{df}

Consider the curve $\Gamma := \{ (x,y) \in \mathbb{C}^2 | f(x) =
g(y)\}$, which can be thought of as the the fibre product of the two
maps $f : \mathbb{C} \ra \mathbb{C}$ and $g : \mathbb{C} \ra
\mathbb{C}$ (or as $(f \times g)^{-1} (\Delta \subset \mathbb{C}
\times \mathbb{C})$, $\Delta$ being the diagonal in $\CC \times \CC$).

Denoting the respective branch loci of $f$ and $g$ by $B_f$ and
$B_g$, we have that the singular locus of $\Gamma$ is contained
in the intersection of the respective inverse images of the
branch loci, i.e., $ Sing (\Gamma) \subset p_x ^{-1} (B_f) \times
p_y ^{-1} (B_g)$. Moreover, if  we denote by $C'$ the normalization
of $\Gamma$, then
  the branch locus $B$ of $\pi :
C'  \ra \CC$ (which has degree $deg(f) \cdot deg(g)$)
is equal to  $B_f \cup B_g$.

We have the monodromy homomorphisms $\mu_f : \pi_1(\CC - B_f) \ra
\mathfrak{S}_n$, $\mu_g : \pi_1(\CC - B_g) \ra
\mathfrak{S}_m$ and via the epimorphism $\pi_1(\CC - B) \ra \pi_1(\CC
- B_f)$ (resp. the one for $g$), we get
$$
\Phi = \mu_f \times \mu_g : \pi_1(\CC - B) \ra \mathfrak{S}_n
\times \mathfrak{S}_m.
$$

Therefore Riemann's existence theorem and the fact that a variety
is irreducible iff its nonsingular locus  is connected yields the following:

\begin{prop}
$f(x) - g(y)$ is irreducible if and only if the product of the two
monodromies $\Phi = \mu_f \times \mu_g$ is transitive.
\end{prop}

An easier irreducibility criterion is obtained by the following observation.

Homogenizing the equation $f(x) = g(y)$ as $ y_0^m x_0^n f(\frac{x_1}{x_0}) =
  x_0^n  y_0^m g(\frac{y_1}{y_0})$ we obtain a compactification
$\overline {\Gamma} \subset \PP^1 \times \PP^1$ such that
the only point at infinity is $(\infty, \infty ) = ((0,1)(0,1))$.

In this point there are local holomorphic coordinates $(u,v)$
such that the local equation of $\overline {\Gamma}$
reads out as $ u^n = v^m$. Set $ d : = G.C.D. (n,m)$:
then  we have  $ u^n - v^m = \prod_{i=1}^d ( u^{\frac nd} - \zeta^i 
v^{\frac md})$,
$\zeta$ being as usual a primitive $d$-th root of $1$.
Thus we obtain

\begin{prop}\label{diffpol}
The curve $\overline {\Gamma}$ has exactly $d$ branches at
infinity, and is in particular irreducible if $ d=1$.
Letting $C$ be the normalization of $\overline {\Gamma}$,
the algebraic function
$f(x) = g(y)$ has branch locus equal to $B_f \cup B_g \cup \{ \infty\}$
if $ n, m \geq 1$.
\end{prop}

A trivial example of a non irreducible difference polynomial is given
by $f=g$, since then $(x-y) | (f(x) - f(y))$.

A less trivial
counterexample is given by the Chebycheff polynomial $T_{2n}(x) +
T_{2n}(y)$ (cf. \cite{davenport}). In fact, $T_{2n}(x) +
T_{2n}(y) = T_{2n}(x) -
( - T_{2n}(y))$ and the monodromy of $( - T_{2n}(y))$ is simply 
obtained by exchanging
the roles of $\tau_1, \tau_2$.

It is easier to view the dihedral monodromies of $T_{2n}(x)$,
  respectively $( - T_{2n}(y))$, as acting on $ \ZZ / 2n$:
$$ \tau_1 (i, j) =( - i, -j -1),\ \  \tau_2 (i, j) = (- i  - 1, -j).$$

We see immediately that  $ \tau_1  \tau_2 (i, j) =(  i + 1, j + 1)$, thus
the cyclic subgroup of $ D_{2n}$ operates trivially on $ (i-j)$
and we conclude that we have exactly $n$ orbits of
cardinality
$ 4n$ in the product $ ( \ZZ / 2n)^2$, and correspondingly a 
factorization of $T_{2n}(x) + T_{2n}(y)$
in $n$ irreducible factors.

Schinzel asked (cf.  \cite{schinzel}) whether one could give examples
without using the arithmetic of Chebycheff polynomials.
Riemann's existence provides plenty of examples, as the following
one where two degree $7$ polynomials yield a difference polynomial
with exactly two irreducible factors, of respective bidegrees $(4,4), (3,3)$.

\begin{prop}
Consider the following three matrices in $ GL (3, \ZZ/2)$,

$
A_1=\begin{pmatrix}
1& 1 & 0 \\
0& 1 & 0  \\
0 & 0 & 1\\
\end{pmatrix},
$
$
A_2=\begin{pmatrix}
1& 0 & 0 \\
1& 1 & 0  \\
0 & 0 & 1\\
\end{pmatrix},
$
$
A_3=\begin{pmatrix}
1& 0 & 0 \\
0& 0 & 1  \\
0 & 1 & 0\\
\end{pmatrix},
$
and let $f$ be the degree $7$ polynomial associated to
the monodromy action of $\pi_1 ( \CC \setminus \{ 0,1, \lambda \})$
determined by
the above three matrices
on the Fano projective plane associated to the vector space  $ V : = 
( \ZZ/2)^3$ .

Let moreover $g$ be the degree $7$ polynomial associated to
the monodromy action of $\pi_1 ( \CC \setminus \{ 0,1, \lambda \})$
on the dual projective plane (associated to the vector space  $ V^{\vee} :
= ( \ZZ/2)^3$) determined by the (inverses of
the) transposes  of the above  three matrices, which have order $2$.

Then the difference polynomial $ f(x) - g(y)$ is the product of two
  irreducible factors, of respective bidegrees $(4,4), (3,3)$.
\end{prop}

\Proof
Each linear map $A_i$ has a fixed subspace of dimension $2$,
and has order two: hence the associated    permutation is a double 
transposition.

One sees easily that the $A_i$'s generate a transitive subgroup,
therefore the associated covering yields the class of a polynomial
$f$ of degree $7$, and similarly for $g$. By construction,
the product action on $ \PP (V) \times \PP (V^{\vee})$
leaves invariant the incidence correspondence, which is
a correspondence of type $(3,3)$.

Moreover, one sees immediately that $A_1, A_2 = ^t A_1$  leave the
vector $e_3$ fixed, and they
operate transitively on the projective line generated by $ e_1, e_2$.
Therefore the action on the product  $ \PP (V) \times \PP (V^{\vee})$
is easily seen to have exactly two orbits, the incidence correspondence
and its complementary set.

\QED

Another very interesting occurrence of Chebycheff  polynomials
concerns {\em Schur's problem}, which asks: given a polynomial $f \in \CC[z]$,
when  is $\frac{f(x) -
   f(y)}{x-y}$  irreducible?

Nontrivial non irreducible Schur polynomials are obtained again through the
Chebycheff polynomials: $\frac{T_n(x) -
   T_n(y)}{x-y}$. Let us briefly explain how,
  using again the description of the product monodromy
(this procedure shows how one can find many more examples
  of non irreducible Schur
polynomials, although without an explicit determination of the coefficients
of the polynomial $f$).

Using the previous notation, we get this time the action
$$ \tau_1 (i, j) =( - i, -j ),\ \   \tau_2 (i, j) = ( -i  - 1, -j-1) 
\ on \ (\ZZ / n)^2,$$
thus in particular $ \tau_1 \tau_2 (i, j) = (i + 1, j + 1)$
and obviously the difference $i-j$ is only transformed into
$ \pm (i-j)$.

Thus the  corresponding Schur polynomial has exactly  $\frac {n-1}2$ 
factors for
$n$ odd, and $\frac {n}2$ factors for
$n$ even.

We already said that Chebycheff polynomials are quite
ubiquitous.

Let us point out another beautiful application of
Chebycheff polynomials: namely,
Chmutov (cf. \cite{chmutov}) used them to construct surfaces in 
$\PP^3$ with ``many''
nodes, obtaining the best known asymptotic lower bounds.

\section{Belyi's Theorem and Triangle curves}

Grothendieck's enthusiasm was raised by the following result, where
Belyi made a very clever
and very simple use of the Belyi polynomials in order to reduce the number of
critical values of an algebraic function

\begin{teo} (Belyi, cf. \cite{belyi})
An algebraic curve $C$ can be defined over $\overline{\QQ}$ if and only
if there exists a holomorphic map $f : C \ra \PP^1_{\CC}$ branched
exactly in $\{0, 1, \infty \}$.
\end{teo}

Again here the monodromy of $f$ is determined by the children's drawing
$ f^{-1} ([0,1])$.
If, moreover, we assume that $f$ is Galois, then we call $C$ a
triangle curve.

\begin{df}
$C$ is a {\em triangle curve} if there is a finite group $G$ acting effectively
on $C$ and with the property that $ C / G \cong \PP^1_{\CC}$, and
$ f : C \ra \PP^1_{\CC} \cong C / G $ has $\{0, 1, \infty \}$ as branch set.
\end{df}

Belyi's  construction of triangle curves is rather complicated. Easier examples
can be constructed using, as in the previous section, difference polynomials
associated (as in \ref{diffpol}) to polynomials with $0,1$ as only 
critical values.

Gabino Gonzalez (cf. \cite{gonz}) was recently able to extend Belyi's 
theorem to the case of
complex surfaces (in terms of Lefschetz maps with three critical values).

In the next section we shall describe some higher dimensional analogues
of triangle curves.

\section{Beauville surfaces}

Inspired by a construction of A. Beauville (cf \cite{bea}) of a 
surface with $K^2 =
8$, $p_g = q=0$ as a quotient of the product of two Fermat curves of
degree $5$ by the action of $\ZZ / 5 \ZZ$ the second author gave in
  \cite{cat00} the following definition.

\begin{df}
A {\em Beauville \ surface} is a compact complex surface $S$ which

1) is {\em rigid}, i.e., it has no nontrivial deformation,

2) is {\em isogenous to a (higher) product}, i.e., it is a quotient
$S = (C_1 \times C_2) /G$ of a product of curves of resp. genera
$\geq 2$ by the
free action of a finite group $G$.
\end{df}

In the above cited paper, the second author introduced and studied
extensively surfaces isogenous to a higher product. Among others it is
shown there that the topology of a surface isogenous to a product
determines its deformation class up to complex conjugation. The
following theorem contains a correction to thm. 4.14 of \cite{cat00}.

\begin{teo}
Let $S = (C_1 \times C_2) /G$ be a surface isogenous to a product.
Then any surface $X$ with the
same topological Euler number and the same fundamental group as $S$
is diffeomorphic to $S$. The corresponding subset of the moduli space
$\mathcal{M}^{top}_S = \mathcal{M}^{diff}_S$, corresponding to 
surfaces homeomorhphic,
resp, diffeomorphic to $S$, is either irreducible and connected
or it contains
two connected components which are exchanged by complex
conjugation.

In particular, if $X$ is orientedly diffeomorphic to $S$, then $X$ is
deformation equivalent to $S$ or to $\bar{S}$.
\end{teo}

The class of surfaces isogenous to a product and their higher
dimensional analogues provide a wide class of examples where one can
test or disprove several conjectures and questions
(cf. e.g. \cite{FabReal}, \cite{torelli}, \cite{beauvillesurf}).

Notice, that given a surface $S = (C_1 \times C_2) /G$ isogenous to a
product, we obtain always three more, exchanging $C_1$ with its
conjugate curve $\bar{C}_1$, or $C_2$ with $\bar{C}_2$, but only if we
conjugate both $C_1$ and $C_2$, we obtain an orientedly diffeomorphic
surface. However, these four surfaces could be all biholomorphic to
each other.

If $S$ is a Beauville surface (and $X$ is orientedly diffeomorphic to
$S$) this implies: $X \cong S$ or $X \cong \bar{S}$.
In other words, the corresponding subset of the
moduli space $\mathcal{M}_S$ consists of one or two points (if we
insist on keeping the orientation fixed, else we may get up to four points).

The interest for Beauville surfaces comes from the fact that they are
the rigid ones amongst surfaces isogenous to a product.
We recall that an algebraic variety $X$ is {\em rigid} if and only if
it does not have any non trivial deformation (for instance, the 
projective space is
rigid). There is another (stronger) notion of rigidity, which is the following

\begin{df}
An algebraic variety $X$ is called {\em strongly rigid} if any other
   variety homotopically equivalent to $X$ is either biholomorphic or
antibiholomorphic to $X$.
\end{df}

\begin{oss}
1) It is nowadays wellknown that smooth compact quotients of symmetric
   spaces are rigid (cf. \cite{calabivesentini}).

2) Mostow (cf. \cite{mostow} proved that indeed locally symmetric
spaces of complex dimension $\geq 2$ are strongly rigid, in the sense
that any homotopy equivalence is induced by a unique isometry.
\end{oss}

These varieties are of general type and the moduli space of varieties
   of general type is defined over $\mathbb{Z}$, and naturally the
absolute Galois group $Gal(\overline{\mathbb{Q}} / \mathbb{Q})$ acts on
the set of their connected components. So, in our special case,
$Gal(\overline{\mathbb{Q}} / \mathbb{Q})$ acts on the isolated points
which parametrize rigid varieties.

In particular, rigid varieties are defined over a number field
and work of Shimura gives in special cases a way of computing
explicitly their fields of definition. By this reason these varieties
were named {\em Shimura varieties} (cf. Deligne's Bourbaki seminar
\cite{delignebourb}).

A quite general question is

\begin{question}
What are the fields of definition of rigid varieties? What is the
$Gal(\overline{\mathbb{Q}} / \mathbb{Q})$-orbit of the point in the moduli
space corresponding to a rigid variety?
\end{question}

Coming back to Beauville surfaces we observe the following:

\begin{oss}
The rigidity of a Beauville surface is equivalent to the  condition that
$(C_i, G^0)$ is a triangle curve, for $i=1,2$ ($G^0 \subset G$ is the subgroup
of index $\leq 2$ which does not exchange the two factors).

It follows that a Beauville surface is defined over $\overline{\QQ}$,
whence the Galois group $Gal (\overline{\QQ} / \QQ)$ operates on the
discrete subset of the moduli space $\mathcal{M}_S$ corresponding
to Beauville surfaces.

By the previous theorem, the Galois group $Gal (\overline{\QQ}/ \QQ)$
may transform a Beauville surface into another one
with a non isomorphic fundamental group.
Phenomena of this kind were already observed by J.P. Serre (cf. \cite{serre}).
\end{oss}

It looks therefore interesting to
   address the following problems:

\begin{question}
Existence and classification of Beauville surfaces, i.e.,

a) which finite groups $G$ can occur?

b) classify all possible Beauville surfaces for a given finite group
$G$.
\end{question}

\begin{question}
Is the Beauville surface $S$ biholomorphic
   to its complex conjugate surface $\bar{S}$?

Is $S$ real (i.e., does there exist a biholomorphic map $\sigma : S
\rightarrow \bar{S}$ with $\sigma^2 = id$)?
\end{question}

The major motivation to find these surfaces is rooted in the following
so called {\em Friedman-Morgan speculation (1987)}
(cf. \cite{friedmor}), which we briefly explain in the following.

One of the fundamental problems in the theory of complex algebraic
surfaces is to understand the moduli spaces of surfaces of general
type, and in particular their connected components, which parametrize
the deformation equivalence classes of minimal surfaces of general
type.

\begin{df}
Two minimal surfaces $S$ and $S'$ are said to be {\em def -
    equivalent} (we also write: $S \sim_{def} S'$) if and only if they
are elements of the same connected component of the moduli space.
\end{df}

By the classical theorem of Ehresmann (\cite{ehre}), two def - 
equivalent algebraic
surfaces are (orientedly) diffeomorphic.

In the late eighties Friedman and Morgan (cf. \cite{friedmor}) conjectured
that two algebraic surfaces are diffeomorphic if and only if they are
def - equivalent ({\bf def = diff}).

Donaldson's breaktrough results had made clear that diffeomorphism and
homeomorphism differ drastically for algebraic surfaces
(cf. \cite{donaldson}) and the success of gauge theory led Friedman and
Morgan to ``speculate'' that the diffeomorphism type of algebraic
surfaces determines the deformation class. After the first
counterexamples of M. Manetti (cf. \cite{man4}) appeared, there
were further counterexamples given by Catanese, Kharlamov-Kulikov,
Catanese-Wajnryb, Bauer-Catanese-Grunewald (cf. \cite{FabReal},
   \cite{k-k},\cite{beauvillesurf}, \cite{catwaj}).

The counterexamples are of quite different nature:  Manetti used $(\Z 
/2)^r$-covers of  $\PP ^1_{\C} \times \PP ^1_{\C}$, his
surfaces have $b_1 =0$, but are not $1$-connected.
Kharlamov - Kulikov used quotients $S$ of the unit ball in $\C^2$, 
thus $b_1 > 0$, and
the surfaces are rigid.
Catanese used surfaces isogenous to a product (which are not rigid).

Our examples are given by Beauville surfaces. The great advantage of
these is that, being isogenous to a product, they can be described by
combinatorial data and the geometry of these surfaces is encoded in
the algebraic data of a finite group.

In order to reduce the description of Beauville surfaces to some
group theoretic statement, we need to recall that surfaces
isogenous to a higher product belong to two types:

\begin{itemize}
\item
$S$ is of {\em unmixed type} if the action of $G$ does not mix the two factors,
i.e., it is the product action of respective actions of $G$ on $C_1$,
resp. $C_2$.
\item
$S$ is of {\em mixed type}, i.e., $C_1$ is isomorphic to $C_2$, and the
subgroup $G^0$ of transformations in $G$ which do not mix the factors
has index precisely $2$ in $G$.
\end{itemize}

The datum of a Beauville surface can be completely described group
theoretically,
since it is equivalent to the datum of two triangle curves with
isomorphic groups with some additional condition assuring that the
diagonal action on the product of the two curves is free.

\begin{df}
Let $G$ be a finite group.

1) A quadruple $v=(a_1,c_1;a_2,c_2)$
of elements of $G$ is an {\em unmixed Beauville structure for $G$}
if and only if

(i) the pairs $a_1,c_1$, and $a_2,c_2$ both generate $G$,

(ii) $\Sigma(a_1,c_1) \cap \Sigma(a_2,c_2) =\{ 1_{G}\}$, where
$$
\Sigma(a,c):=\bigcup_{g\in G}
\bigcup_{i=0}^\infty\ \{ga^ig^{-1},gc^ig^{-1},g(ac)^ig^{-1}\}.
$$
We write $\mathbb{B}(G)$ for the set of unmixed Beauville structures on $G$.

2) {\em A mixed Beauville
quadruple for $G$} is a quadruple
$M=(G^0;a,c;g)$ consisting of a
subgroup $G^0$ of index $2$ in $G$, of elements
$a,c\in G^0$ and of an element $g\in G$ such that

i) $G^0$ is generated by $a,c$,

ii) $g\notin G^0$,

iii) for  every $\gamma \in G^0$ we have
$g\gamma g\gamma\notin \Sigma(a,c)$.

iv) $\Sigma(a,c)\cap\Sigma(gag^{-1},gcg^{-1})=\{ 1_G \}$.

\noindent
We call $\M (G)$ the set of mixed Beauville quadruples on the
group $G$.
\end{df}

\begin{oss}
We consider here finite groups $G$ having a pair $(a,c)$ of generators. Setting
$(r,s,t):= ({\rm ord}(a),{\rm ord}(c),{\rm ord}(ac))$, such a group is a
quotient of the triangle group
$$
T(r,s,t):=\langle x,y \ |\ x^r=y^s=(xy)^t=1\rangle.
$$
\end{oss}

We defined some algebraic structures on a finite group in the above. As
the name ``Beauville structure'' already suggests, these data give
rise to a Beauville surface as follows:

Take as base point
$\infty \in \PP^1$, let $B : = \{-1, 0, 1\}$, and  take
the following generators $\alpha, \beta$ of
$\pi_1 (\PP^1_{\C} - B, \infty)$ ($\gamma := (\alpha \cdot \beta )^{-1}$):

\font\mate=cmmi8
\def\min{\hbox{\mate \char60}\hskip1pt}

\begin{figure}[htbp]
\begin{center}
\scalebox{1}{\includegraphics{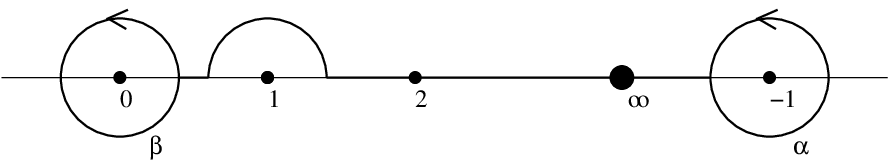}}
\end{center}
\label{figura9}
\end{figure}

We observe that we now prefer to take as branch locus the set $\{-1,
0, 1\}$, since it is more convenient for describing the complex
conjugation group theoretically.

Let now $G$ be a finite group and $v=(a_1,c_1;a_2,c_2) \in \mathbb{B}(G)$.
We get surjective homomorphisms
\begin{equation}\label{e1}
\pi_1 (\PP^1_{\C} - B, \infty)\to G,\qquad \alpha\mapsto a_i, \
\gamma \mapsto c_i
\end{equation}
and Galois coverings
$\lambda_i:C(a_i, c_i) \rightarrow \PP^1_\C$
ramified only in $\{-1,0,1\}$ with ramification indices equal to the
orders of $a_i$, $b_i$, $c_i$ and with group $G$) ({\em Riemann's
existence theorem}).

\begin{oss}
1) Condition (1), ii) assures that the action of $G$ on $C(a_1, c_1)
\times C(a_2, c_2)$ is free.

2) Let be $\iota(a_1,c_1;a_2,c_2) =
(a_1^{-1},c_1^{-1};a_2^{-1},c_2^{-1})$. Then $S(\iota(v)) =
\overline{S(v)}$. (Note that $\bar{\alpha} = \alpha^{-1}$,
$\bar{\gamma} = \gamma^{-1}$.)

3) We have: $g(C(a_1,c_1))\ge 2$ and
$g(C(a_2,c_2))\ge 2$. This is a nontrivial fact, whcih comes from
group theory. It is equivalent to the fact that $\mu(a_i,c_i):=
\frac{1}{ord(a_i)}+\frac{1}{ord(c_i)}+\frac{1}{ord(a_ic_i)}<1$
(cf. prop. \ref{hyperbolic}).
\end{oss}

In order to address reality questions of Beauville surfaces, we have to
translate the two questions:
\begin{itemize}
\item Is $S$ biholomorphic to its complex conjugate $\bar{S}$?
\item Is $S$ real, i.e., does there exist a biholomorphism $\sigma : S
   \ra \bar{S}$, such that $\sigma^2 = id$?
\end{itemize}
into  group theoretic conditions.

\begin{oss}
Actually we can define
a finite permutation group ${\rm A}_\mathbb{B}(G)$ such that for
$v,v'\in\mathbb{B}(G)$ we have : $S(v) \cong S(v')$ if and only if $v$ is
in the
${\rm A}_\mathbb{B}(G)$-orbit of $v'$. Under certain assumptions on
the orders of the generating elements these orbits are easy to
describe as shows the following result.
\end{oss}

\begin{prop}
Let $G$ be a finite group and $v=(a_1,c_1;a_2,c_2) \in\mathbb{B}(G).$

Assume that $\{ {\rm ord}(a_1), {\rm ord}(c_1),{\rm ord}(a_1c_1)\}
\neq \{ {\rm ord}(a_2), {\rm ord}(c_2),{\rm ord}(a_2c_2)\}$ and that
${\rm ord}(a_i) < {\rm ord}(a_i c_i) < {\rm ord}(c_i)$.

Then $S(v) \cong \overline{S(v)}$
if and only if there are  inner automorphisms $\phi_1$, $\phi_2$ of
$G$ and an automorphism $\psi \in Aut(G)$
such that, setting
$\psi_j : = \psi \circ \phi_j$, we have
$$
\psi_1(a_1) = a_1^{-1}, \ \ \ \psi_1 (c_1) = c_1^{-1},
$$
$$
\psi_2 (a_2) = {a_2}^{-1}, \ \ \ \psi_2 (c_2) = {c_2}^{-1}.
$$
In particular  $S(v)$ is isomorphic to $\overline{S(v)}$ if and only if
$S(v)$ has a real structure.
\end{prop}

The following is an immediate consequence of the above
\begin{oss}
If $G$ is abelian, $v \in\mathbb{B}(G)$. Then $S(v)$ always has a real
structure, since $g \mapsto -g$ gives an automorphism of $G$ of order two.
\end{oss}

Concerning the existence of (unmixed) Beauville groups respectively
unmixed Beauville surfaces we have among others the following results:

\begin{teo}
1) A finite abelian group $G$ admits an unmixed Beauville structure iff $G
\cong (\mathbb{Z} /n )^2$, $(n,6) = 1$.

\noindent
2) The following groups admit unmixed Beauville structures:

a) $\mathfrak{A}_n$ for large $n$,

b) $\mathfrak{S}_n$ for $n\in \NN$ with $n\ge 7$,

c) ${\bf SL}(2,\F_p)$, ${\bf PSL}(2,\F_p)$ for $p \neq 2,3,5$.
\end{teo}

We checked all finite simple nonabelian groups of order $\leq 50000$
and found unmixed Beaville structures on all of them with the
exception of $\mathfrak{A}_5$. This led us to the following
\begin{conj} All finite simple nonabelian groups except $\AAA_5$
admit an unmixed Beauville structure.
\end{conj}

We have checked this conjecture for some bigger simple groups like the
Mathieu groups ${\bf M}12, \,{\bf M}22 $ and also matrix groups of
size bigger than $2$.

We call $(r,s,t)\in \NN^3$ {\it hyperbolic} if
$$\frac{1}{r}+\frac{1}{s}+\frac{1}{t} < 1.$$
In this case the triangle group $T(r,s,t)$ is hyperbolic.
   From our studies also the following looks suggestive:

\begin{conj}
Let  $(r,s,t)$,  $(r',s',t')$ be two hyperbolic
     types. Then almost all alternating groups $\AAA_n$ have an 
unmixed Beauville
     structure $v=(a_1,c_1;a_2,c_2)$ where $(a_1,c_1)$ has type  $(r,s,t)$ and
    $(a_2,c_2)$ has type  $(r',s',t')$.
\end{conj}

The above conjectures are variations of a conjecture of Higman
(proved by B. Everitt (2000), \cite{everitt}) asserting that
every hyperbolic triangle group surjects onto almost all alternating
groups.

The next result gives explicit examples of rigid surfaces not
biholomorphic to their complex conjugate surface.
\begin{teo}
The following groups admit unmixed Beauville structures $v$ such that $S(v)$ is
not biholomorpic to $\overline{S(v)}$:

\begin{itemize}
\item[1)] the symmetric group $\mathfrak{S}_n$ for $n\ge 7$,
\item[2)] the alternating group $\mathfrak{A}_n$ for $n\ge 16$ and 
$n\equiv 0$ mod $4$,
$n\equiv 1$ mod $3$, $n\not\equiv 3,4$ mod $7$.
\end{itemize}
\end{teo}

The following theorem gives examples of surfaces which are not real,
but biholomorphic to their complex conjugates, or in other words, they
give real points in their moduli space which do not correspond to 
real surfaces.

\begin{teo}\label{nonreal}
Let $p>5$ be a prime with $p\equiv 1$ mod $4$, $p\not\equiv 2,4$ mod 5,
$p\not\equiv 5$ mod $13$ and $p\not\equiv 4$ mod $11$. Set $n:=3p+1$.
Then
there is  an unmixed Beauville surface $S$ with group $\mathfrak{A}_n$
which is biholomorphic to the complex conjugate surface $\bar{S}$, but is not
real.
\end{teo}

For {\em mixed} Beauville surfaces the situation is far more complicated,
as already the following suggests.
\begin{teo}
1) If a finite group $G$ admits a mixed Beauville structure, then the subgroup
$G^0$ is non abelian.

2) No group of order $\leq 512$ admits a mixed Beauville structure.
\end{teo}

We give a general construction of finite groups admitting a mixed
Beauville structure.


Let $H$ be a non-trivial group. Let
$\Theta :H\times H \to H\times H$ be the automorphism defined by
$\Theta(g,h):=(h,g)$ ($g,h\in H$). We consider the semidirect product
\begin{equation}
H_{[4]}:= (H\times H) \rtimes  \ZZ / 4\ZZ
\end{equation}
where the generator $1$ of $ \ZZ / 4\ZZ$ acts through $\Theta$ on $H\times H$.
Since $\Theta^2$ is the identity we find
\begin{equation}
H_{[2]}:=H\times H\times 2\ZZ / 4\ZZ \cong H\times H\times \ZZ / 2\ZZ
\end{equation}
as a subgroup of index $2$ in $H_{[4]}$.

We have now
\begin{lem}
Let $H$ be a non-trivial finite group and let $a_1,c_1$,  $a_2,c_2$ 
be elements of
$H$. Assume that

1. the orders of $a_1,c_1$ are even,

2. $a_1^2, a_1c_1, c_1^2$ generate $H$,

3. $a_2,c_2$ also generate $H$,

4.$({\rm ord}(a_1) \cdot {\rm ord}(c_1) \cdot {\rm ord}(a_1c_1),{\rm
ord}(a_2) \cdot {\rm ord}(c_2) \cdot {\rm ord}(a_2c_2)) = 1$.

Set $G:=H_{[4]}$, $G^0:=H_{[2]}$ as above
and $a:=(a_1,a_2,2)$, $c:=(c_1,c_2,2)$. Then
$(G^0;a,c)$ is a mixed Beauville structure on $G$.
\end{lem}

\Proof It is easy to see that  $a,c$ generate  $G^0:=H_{[2]}$.

The crucial observation is
\begin{equation}
(1_H,1_H,2)\notin \Sigma(a,c).
\end{equation}
If this were not correct, it would have to be  conjugate of a power
of $a$, $c$ or $b$. Since the
orders of $a_1$, $b_1$, $c_1$ are even, we obtain a contradiction.

Suppose that
$h=(x,y,z)\in\Sigma(a,c)$ satisfies
${\rm ord}(x)={\rm ord}(y)$: then our condition 4 implies that $x=y=1_H$ and
(5) shows $h=1_{H_{[4]}}$.

Let now $g\in H_{[4]}$, $g\notin H_{[2]}$
and $\gamma\in G^0=H_{[2]}$ be given. Then $g\gamma=(x,y,\pm 1)$ for
appropriate $x,y\in H$. We find
$$ (g\gamma)^2=(xy,yx,2)$$
and the orders of the first two components of $(g\gamma)^2$ are the
same, contradicting the above remark.

Therefore the third condition is satisfied.

We come now to the fourth condition of a
mixed Beauville quadruple. Let $g\in H_{[4]}$, $g\notin H_{[2]}$ be given,
for instance $(1_H,1_H,1)$.
Conjugation with $g$ interchanges then
the first two components of an element $h\in H_{[4]}$.
Our hypothesis 4 implies the result. \QED

As an application we find the following examples
\begin{teo}
    Let $p$ be a prime with $p\equiv 3$ mod $4$ and
    $p\equiv 1$ mod $5$ and consider the group $H:={\bf SL}(2,\F_p)$.
Then $H_{[4]}$ admits a mixed Beauville structure $u$ such
that $S(u)$ is not biholomorphic to $\overline{S(u)}$.
\end{teo}

\begin{oss}
Note that the smallest prime satifying the
above congruences is $p = 11$ and we get that $G$ has order equal to $6969600$.
\end{oss}

So it is natural to ask the following:

\begin{question}
What is the minimal order of a group $G$ admitting a mixed Beauville structure?
\end{question}

It is interesting to observe that one important numerical restriction
follows automatically from group theory:

\begin{prop}\label{hyperbolic}
Assume $(a_1,c_1;a_2,c_2)\in\mathbb{B}(G)$. Then $\mu(a_1,c_1) :=
\frac{1}{{\rm ord}(a_1)} + \frac{1}{{\rm ord}(c_1)}
+ \frac{1}{{\rm ord}(a_1c_1)}< 1$ and
$\mu(a_2,c_2) < 1$. Whence we have: $g(C(a_1,c_1))\ge 2$ and
$g(C(a_2,c_2))\ge 2$.
\end{prop}

\Proof
We may without loss of generality assume that $G$ is not cyclic.
Suppose $(a_1,c_1)$ satisfies $\mu(a_1,c_1)>1$: then the type of $(a_1,c_1)$ is
up to permutation amongst the
$$(2,2,n)\ (n\in \NN),\ (2,3,3),\ (2,3,4),\ (2,3,5).$$
These can be excluded easily.
If $\mu(a_1,c_1)=1$ then the type of $(a_1,c_1)$ is
up to permutation amongst the
$$(3,3,3),\ (2,4,4),\ (2,3,6)$$
and $G$ is a finite quotient of one of the wall paper groups and cannot admit
an unmixed Beauville structure.
\QED

We finish sketching the underlying idea for the
examples of  Beauville surfaces not isomorphic to their conjugate surface
obtained from symmetric groups.

\begin{lem}
Let $G$ be the symmetric group $\mathfrak{S}_n$ in $n \geq 7$ letters
and let $p$ be an odd prime. We set $a:=(1,p+2,p+1)(2,p+3)$, $c:=(1,2,
\ldots , p)(p+1, \ldots ,n)$. Then the following holds:

1) there is no
automorphism of $G$ carrying $ a \ra a ^{-1}$, $c \ra c^{-1}$;

2) $\mathfrak{S}_n = <a,c>$.

\end{lem}

\Proof 1) Since $n \neq 6$, every automorphism of $G$ is an inner one.
If there is a permutation $g$ conjugating $ a $ to $ a ^{-1}$, $c $ to
$ c^{-1}$, $g$ would leave each of the sets $\{1,p+1,p+2\}$, $\{2,p+3
\}$, $\{1,2,\dots ,p\}$, $\{p+1,\dots ,n\}$ invariant. By looking
at their
intersections we conclude that $g$ fixes the elements
$1,2,p+3$. But then $gcg^{-1} \neq
c^{-1}$.

2) Let $G':= <a,c>$. Note that $G'$ contains the transposition $a^3 =
(2,p+3)$ and it contains the $n$-cycle $c \cdot (2,p+3)$. Therefore
$G' = \mathfrak{S}_n$.
\QED

\begin{oss}
1) $ a' : =  \sigma^{-1}$, $c' : = \tau
\sigma^2$, where
$\tau := (1,2)$
and $\sigma: = (1,2, \dots ,
n)$. It is obvious that $\mathfrak{S}_n = <a',c'>$.

2) Let $n \geq 7$ and let $p$ be an odd prime such that $n$ is not congruent
to $0$ or $1 (mod \ p)$ (cf. \ref{prime}), we get by the same arguments as in
\cite{beauvillesurf}, prop. 5.5. that $\Sigma(a,c) \cap
\Sigma(a',c') = \{1\}$.

In fact, one has to observe that conjugation preserves the types
(coming from the cycle decomposition). The types in $\Sigma(a,c)$ are
derived from $(3)$, $(2)$, $(p)$, $(n-p)$ and $(n-p-1)$ (since we assume that
$p$ does neither divide $n$ nor $n-1$), whereas the types in
$\Sigma(a',c')$ are derived from $(n)$, $(n-1)$, or $(\frac{n-1}{2}, 
\frac{n+1}{2})$.
\end{oss}

The existence of an odd prime $p$ as above is assured by the following
lemma.

\begin{lem}\label{prime}
Let $n \geq 5$ be any integer $\neq 6$. Then there is an odd prime
number $p$ such that $n$ is not congruent to $0$ or $1 (mod \ p)$.
\end{lem}

\Proof
If $n$ is odd, then $n-2 \geq 3$ and we choose an odd prime $p$
dividing $n-2$. Then $n \equiv 2(mod \ p)$ and we are done.

Assume now $n$ to be even and let $p$ be a prime dividing $n-3$. If
there exists such a prime $\neq 3$, then we are done, since $n \equiv
3(mod \ p)$. Otherwise, $n-3 = 3^k$, for some $k \geq 2$. Then it is
obvious that $n+3 = 3^k +6$ is not a power of $3$ and we take a prime
$p>3$ dividing $n+3$ and we are done, since $n \equiv -3(mod \ p)$ and
$-3$ is not congruent to $0$ or $1(mod \ p)$:
\QED

\end{document}